\documentclass[psamsfonts,12pt]{amsart}
\usepackage{amssymb}
\usepackage[left=1in,right=1in,top=1in,bottom=1in]{geometry}

\newtheorem{deff}{Definition}[section]
\newtheorem{lemma}[deff]{Lemma}
\newtheorem{theorem}[deff]{Theorem}
\newtheorem{corollary}[deff]{Corollary}

\newtheorem{fact}[deff]{Fact}

\newtheorem{em-example}[deff]{Example}
\newtheorem{em-def}[deff]{Definition}        
\newtheorem{em-remark}[deff]{Remark}         
\newtheorem{em-question}[deff]{Question}

\newenvironment{example}{\begin{em-example} \em }{ \end{em-example}}
\newenvironment{definition}{\begin{em-def} \em  }{ \end{em-def}}
\newenvironment{remark}{\begin{em-remark} \em }{\end{em-remark}}

\newenvironment{claim}{\underline{Claim}: \em}{\medskip}

\def\ker{\mathop{\rm ker}}

\def\T{{\mathbb T}}
\def\Z{{\mathbb Z}}
\def\N{{\mathbb N}}
\def\R{{\mathbb R}}
\def\Q{{\mathbb Q}}
\def\K{{\widehat{\Q}}}

\def\P{{\mathbb P}}
\def\cont{\mathfrak c}
\def\hull#1{\langle{#1}\rangle}
\def\gr#1{\langle{#1}\rangle}

\title[Super-sequences in the arc component]{Super-sequences in the arc component of a compact connected group}

\author[D. Dikranjan]{Dikran Dikranjan}
\address[Dikran Dikranjan]{Universit\`a di Udine, Dipartimento di Matematica e Informatica
\\
 via delle Scienze, 206 - 33100 Udine, Italy}
\email{dikran.dikranjan@dimi.uniud.it}

\author[D. Shakhmatov]{Dmitri Shakhmatov}
\address[Dmitri Shakhmatov]{Graduate School of Science and Engineering,
Division of Mathematics, Physics and Earth Sciences\\
Ehime University, Matsuyama 790-8577, Japan}
\email{dmitri@dpc.ehime-u.ac.jp}

\begin{document}
\begin{abstract} 
  Let $G$ be an abelian topological group. The symbol $\widehat{G}$
  denotes the group of all continuous characters $\chi:G\to\T$ endowed
  with the compact open topology.  A subset $E$ of $G$ is said to be
  {\em qc-dense in $G$\/} provided that
  $\chi(E)\subseteq\varphi([-1/4,1/4])$ holds only for the trivial
  character $\chi\in \widehat{G}$, where $\varphi: \R \to \T=\R/\Z$ is
  the canonical homomorphism.  A {\em super-sequence\/} is a non-empty
  compact Hausdorff space $S$ with at most one non-isolated point (to
  which $S$ {\em converges\/}).  We prove that an infinite compact
  abelian group $G$ is connected if and only if its arc component
  $G_a$ contains a super-sequence converging to $0$ that is qc-dense
  in $G$.  This gives as a corollary a recent theorem of Au\ss
  enhofer: For a connected locally compact abelian group $G$, the
  restriction homomorphism $r:\widehat{G}\to \widehat{G}_a$ defined by
  $r(\chi)=\chi\restriction_{G_a}$ for $\chi\in\widehat{G}$, is a
  topological isomorphism. We also show that an infinite compact group
  $G$ is connected if and only if its arc component $G_a$ contains a
  super-sequence $S$ converging to the identity $e$ that generates a
  dense subgroup of $G$ (equivalently, $S\setminus\{e\}$ is an
  infinite suitable set for $G$ in the sense of Hofmann and Morris).
\end{abstract}



\dedicatory{Dedicated to Karl H. Hofmann on the occasion of his 76th anniversary}

\thanks{The first author was partially supported by MEC. MTM2006-02036 and FEDER FUNDS}

\thanks{The second 
author was partially supported by the Grant-in-Aid for Scientific Research 
no.~19540092 by the Japan Society for the Promotion of Science (JSPS)}

\maketitle

\section{Introduction} 

{\em All topological groups are assumed to be Hausdorff\/}. Let $G$ be a topological group.
We denote by $\widehat{G}$ the group of all continuous characters $\chi:G\to\T$ endowed with the compact open topology.  A subgroup $D$ of $G$ {\em determines} $G$ if the restriction homomorphism $r:\widehat{G}\to \widehat{D}$ defined by $r(\chi)=\chi\restriction_D$ for $\chi\in\widehat{G}$, is a topological isomorphism \cite{CRT1}. 
If $G$ is locally compact and abelian, then every subgroup $D$ that determines $G$ must be dense in $G$. (When $D$ is dense in $G$, the map $r:\widehat{G}\to \widehat{D}$
is a continuous isomorphism.) The following two theorems are the cornerstone results in the topic of determining subgroups:

\begin{theorem}\label{CA} \cite{Diss,CM}  A metrizable abelian group  $G$ is  determined by each dense subgroup of $G$.  
\end{theorem}

\begin{theorem}
\label{compact:determined:groups:are:metrizable}
{\rm \cite{HMT}}
Every non-metrizable compact  group $G$ contains a dense subgroup that does not 
 determine $G$.\footnote{A proof of this theorem under the assumption of CH can be found in \cite{CRT1}.}
\end{theorem}

According to a well-known  classical result of Eilenberg and Pontryagin,  in a  connected locally compact abelian group $G$ the arc component $G_a$ is dense.
Since a subgroup of a locally compact abelian group determining  it must be dense, the following theorem,  recently proved  by Au\ss enhofer, is a strengthening of this classical result: 

\begin{theorem}\label{MainThm}
\cite{A2}
The arc component $G_a$ of a connected locally compact abelian group $G$ determines $G$.  
\end{theorem}

 While this theorem
is a corollary of Theorem \ref{CA} for a metrizable group $G$, in the non-metrizable case
the mere density of $G_a$ in $G$ ensured by the classical result of Eilenberg and Pontryagin mentioned above, need not guarantee that $G_a$ determines $G$
(as witnessed by Theorem \ref{compact:determined:groups:are:metrizable}).

Theorem \ref{MainThm} is used in \cite{A1} to prove that the uncountable powers of $\Z$ are not strongly reflexive, thereby resolving a problem raised by Banaszczyk on whether uncountable powers of the reals $\R$ are strongly reflexive. 

Let $\varphi: \R \to \T=\R/\Z$ be the canonical homomorphism and $\T_+=\varphi([-1/4,1/4])$.  We will say that a subset $E$ of a topological group $G$ is {\em qc-dense in $G$\/}   (an abbreviation for {\em quasi-convexly dense\/}) provided that $\chi(E)\subseteq \T_+$ only for the trivial continuous homomorphism $\chi: G\to \T$. 

This notion was introduced in \cite{DDL} in the abelian context, and its significance for applications has been recently demonstrated in \cite{DS}.  In particular, qc-density was used in \cite{DS} to establish essential properties of determining subgroups of compact abelian groups, thereby allowing to get a 
short elementary
proof of Theorem \ref{compact:determined:groups:are:metrizable}.

The host of applications of qc-dense sets is made possible by the ultimate connection between the notions of determining subgroup and qc-density described in
the next fact (proved in \cite[Fact 1.4]{DS}). It is a particular case of a more general fact stated without proof (and in equivalent terms) in \cite[Remark 1.2(a)]{CRT1} and \cite[Corollary 2.2]{HMT}.

\begin{fact}  
\label{connection:between:determination:and:qc-density}
A subgroup $D$ of a compact abelian group $G$ determines it if and only if there exists a compact subset of $D$ that is qc-dense in $G$.
\end{fact}

It has been recently shown in \cite{DS} that qc-dense compact subsets (and thus determining subgroups) of a compact abelian group 
must be rather big.

\begin{theorem}
\label{weight:of:qc-dense:set}
\cite[Corollary 2.2]{DS}
If a closed subset $X$ of an infinite compact abelian group $G$ is qc-dense in $G$, then  
$w(X) =  w(G)$. (Here $w(X)$ denotes the weight of a space $X$.)
\end{theorem}

A {\em super-sequence\/} is a non-empty compact Hausdorff space $X$ with at most one non-isolated point $x^*$ \cite{DS1}. 
We will call $x^*$ the {\em limit\/} of $X$ and say that $X$ {\em converges to $x^*$\/}.  Observe that 
a countably infinite super-sequence is a convergent sequence
(together with its limit).

Au\ss{}enhofer \cite{Diss} essentially proved that every infinite compact metric abelian group has a qc-dense sequence converging to 0.\footnote{\label{footnote:no:2}This is an immediate consequence of \cite[Theorem 4.3 or Corollary 4.4]{Diss}. In fact, a more general statement immediately follows from these results: Every dense subgroup $D$ of a compact metric abelian group $G$ contains a sequence converging to 0 that is qc-dense in $G$.}  This result has been recently extended to all compact groups by replacing convergent sequences with super-sequences:  

\begin{theorem}
\cite{DS}
\label{theorem:2}
Every infinite compact abelian group contains a  qc-dense super-sequence converging to $0$.
\end{theorem}

A subspace $X$ of a topological group $G$ {\em topologically generates $G$\/} if the subgroup of $G$  generated by $X$ is dense in $G$. The proof of the following fact is straightforward.

\begin{fact}
\label{qc-dense:implies:topological:generation}
{\rm \cite[Fact 1.3(ii)]{DS}}
Every qc-dense subset of a compact abelian group topologically generates it.
\end{fact}

Let $G$ be a topological group with the identity $e$. If a discrete subset $S$ of $G$ topologically generates $G$ and $S\cup \{e\}$ is closed in $G$, then $S$ is called a {\em suitable set for $G$\/} 
\cite{HM}.

\begin{remark}\label{suitable:remark}
 Clearly, if $S$ is a super-sequence in $G$ that converges to $e$ and topologically generates $G$, then $S\setminus \{e\}$  is a suitable set for $G$. Conversely, if $G$ is compact and $S$ is a suitable set for $G$, then $S\cup \{e\}$ must be a super-sequence. 
\end{remark}

It follows from this remark that a subgroup $D$ of a compact group $G$  contains a super-sequence converging to the identity that topologically generates $G$ if and only if $D$ contains an infinite suitable set for $G$.

Hofmann and Morris discovered the following fundamental result:

\begin{theorem}
\label{Hofmann:Morris}
{\rm (\cite{HM}; see also \cite{HMbook})}
Every compact group has a suitable set.
\footnote{A ``purely topological'' proof of this result based on Michael's selection theorem can be found in \cite{Sh}.}
\end{theorem}  

\begin{remark}\label{Remark:HM}
\begin{itemize}
\item[(i)]
{\em Theorem \ref{theorem:2} implies the particular case of Theorem \ref{Hofmann:Morris}
for abelian groups.\/} Indeed,
let $G$ be a compact abelian group. If $G$ is finite, then $G$ is discrete, and so $G$ itself is a suitable set for $G$. 
Assume now that $G$ is infinite. By Theorem \ref{theorem:2}, $G$ contains a qc-dense super-sequence $S$ converging to $0$. By Fact \ref{qc-dense:implies:topological:generation}, $S$ topologically generates $G$. According to Remark \ref{suitable:remark}, $S\setminus\{0\}$ is a suitable set for $G$.

\item[(ii)] {\em A suitable set for a compact abelian group $G$ need not be qc-dense in $G$\/}.
Indeed, it is well-known that the group $\T^\cont$ is monothetic, that is, topologically generated by a singleton $S$. 
Clearly, $S$ is a suitable set for $\T^\cont$.
Since $w(S)\le\omega<\cont=w(\T^\cont)$, $S$ cannot be qc-dense in $\T^\cont$  by Theorem \ref{weight:of:qc-dense:set}. 

\item[(iii)] It follows from item (ii) that {\em the particular case of Theorem \ref{Hofmann:Morris}
for abelian groups does not imply Theorem \ref{theorem:2}\/}.

\end{itemize}
\end{remark}

\section{Results}

Our first result is  a particular version of Theorem \ref{theorem:2} that characterizes connected compact abelian groups.

\begin{theorem}
\label{super-sequence:theorem}
For  an infinite compact abelian group $G$ the following conditions are equivalent:
\begin{itemize}
\item[(i)]
the arc component $G_a$  of $G$ contains a super-sequence  converging to $0$ that is qc-dense in $G$;
\item[(ii)]  $G$ is connected.
\end{itemize}
\end{theorem}

In view of Fact \ref{connection:between:determination:and:qc-density},  the implication (ii)$\to$(i) of Theorem \ref{super-sequence:theorem}  yields Theorem \ref{MainThm} when $G$ is compact. The general case of Theorem \ref{MainThm} easily follows from the compact case, see the  proof in the end of Section \ref{proofs}.  Therefore, as a by-product, our proof of Theorem \ref{super-sequence:theorem} also provides an  alternative short and self-contained proof of the theorem of Au\ss enhofer.

\begin{remark}
 In view of Theorems \ref{theorem:2} and \ref{super-sequence:theorem}, as well as Au\ss{}enhofer's result cited in footnote \ref{footnote:no:2}, the reader may wonder if every dense subgroup of a compact abelian group $G$ contains a super-sequence converging to 0 that is qc-dense in $G$.  The answer to this question is negative: 
{\em Every non-metrizable compact abelian group $G$ contains a dense subgroup $H$ such that no super-sequence $S\subseteq H$ is qc-dense in $G$\/}. Indeed, apply Theorem \ref{compact:determined:groups:are:metrizable} to get  a dense subgroup $H$ of $G$ that does not determine $G$.
If $S\subseteq H$ is a super-sequence, then (being compact) $S$ cannot be qc-dense in $G$ by Fact \ref{connection:between:determination:and:qc-density}.
\end{remark}

Our second result is a particular version of Theorem \ref{Hofmann:Morris} that characterizes connected compact groups:  

\begin{theorem}
\label{topologically:generating:sequence}
For an infinite compact group $G$ the following conditions are equivalent:
\begin{itemize}
\item[(i)] the arc component $G_a$ of $G$ contains a suitable set for $G$;
\item[(ii)] $G_a$ contains an infinite suitable set for $G$ that is qc-dense in $G$;
\item[(iii)] $G$ is connected.
\end{itemize}
\end{theorem}

\begin{remark}  
One {\em cannot\/} add the following item to the list of equivalent conditions in  Theorem \ref{topologically:generating:sequence}:
\begin{itemize}
\item[(iv)] $G_a$ contains a super-sequence converging to the identity that is qc-dense in $G$. 
\end{itemize}
Indeed, for every finite simple non-commutative group $L$ the compact group $G= \T \times L$ has a sequence $S\subseteq G_a= \T\times \{e\}$ converging to the identity  $(0,e)$ of $G$ that is qc-dense in $G$, see Example  \ref{Example:2}. Since $G$ is not connected, the implication (iv)$\to$(iii) fails.
\end{remark}

\begin{remark}
Let $\cont$ denote the cardinality of the continuum.
{\em There exists a dense (connected, locally connected, countably compact) subgroup $H$ of (the compact, connected abelian group) $G=\T^{2^\cont}$ such that  $H$ contains no suitable set for $G$\/}. Indeed,  one can take as $H$ the dense subgroup of $G$ without a suitable set  for $H$
constructed in \cite[Corollary 2.9]{DTT}. Assume that $S\subseteq H$ is a  suitable set for $G$.
Then $S$ is discrete and $S\cup\{0\}$ is closed in $G$. Since $S\cup\{0\}\subseteq H$, it follows that $S$ is closed in $H$.  Since $S$ topologically generates $G$, it topologically generates $H$ as well. Hence, $S$ is a suitable set for $H$, a contradiction.
\end{remark}

\section{A qc-dense super-sequence in the arc component of $\widehat{\mathbb{Q}}$}

Our main result in this section is  Lemma \ref{lemma3}.  It follows from the density of $\K_a$ in $\K$ and the general result of Au\ss{}enhofer quoted in the 
footnote \ref{footnote:no:2}. However, Au\ss{}enhofer's proof  relies on Arzela-Ascoli theorem and an inductive construction, so the qc-dense sequence she constructs in her proof is  ``generic''. To keep this manuscript self-contained, we provide  a ``constructive'' example of a ``concrete'' qc-dense sequence in $\mathbb{Q}_a$. 

The proof of the following fact is straightforward from the definition.

\begin{fact}\label{continuous:images:are:qc-dense}
Let  $G$ and $H$ be topological groups and $\pi: H\to G$  a continuous surjective group homomorphism. If  a subset $E$ of $H$ is qc-dense in $H$, then $\pi(E)$ is qc-dense in $G$. 
\end{fact}

In the sequel $\N$ denotes the set of natural numbers. 

\begin{example} \label{Example1} 
Let $T=\left\{\frac{1}{2n}: n\in \N, n\ge 1\right\}\cup\{0\}$. The set {\em $\varphi(T)$ is a qc-dense sequence in $\T$ converging to $0$\/}. 
 Indeed, let $\chi\in\widehat{\T}$ be a non-zero character. Then there exists $m\in\Z\setminus\{0\}$ such that $\chi(x)=mx$ for all $x\in\T$.  Let $n=|m|$. Then $\frac{1}{2n}\in T$ and so $x=\varphi\left(\frac{1}{2n}\right)\in\varphi(T)$. Since $\chi(x)=mx=\varphi\left(\frac{m}{2n}\right)=\varphi\left(\frac{1}{2}\right)\not\in\T_+$, we have $\chi(\varphi(T))\setminus \T_+\not=\emptyset$. This proves that $\varphi(T)$ is qc-dense in $\T$.
\end{example}

For $g\in G$ the symbol $\hull{g}$ denotes the cyclic subgroup of $G$ generated by $g$. 

\begin{lemma}\label{lemma2}
Let $\P=\{p_n:n\in\N\}$ be a faithful enumeration of the set $\mathbb{P}$ of prime numbers. Define 
$$
H=\prod_{n\in \N} \Z_{p_n},
$$
 and let $v=\{1_{p_n}\}_{n\in\N}\in H$, where each $1_{p_n}$ is the identity of $\Z_{p_n}$.  For $n\in\N$ define
$
k_n=(p_0p_1\ldots p_{n-1})^n.
$
Then the set
\begin{equation}
\label{definition:of:S}
S=\{m k_n v: n\in \N, m\leq k_{n+1}\}\cup\{0\}\subseteq \gr{v}
\end{equation}
 is a sequence converging to $0$ that is  qc-dense in $H$.
\end{lemma}

\begin{proof} For $n\in \N$ define
\begin{equation}
\label{definition:of:Wn}
W_n= k_nH= p_0^{n} \Z_{p_0}\times p_1^{n} \Z_{p_1}\times \ldots \times p_{n-1}^{n} \Z_{p_{n-1}}\times \prod_{i=n}^\infty \Z_{p_i}.
\end{equation}
(Note that $k_0=1$.) Then $\{W_n:n\in\N\}$ forms a base of $H$ at $0$ consisting of clopen subgroups. 
It is easy to see that each $W_n$ may miss only finitely many members of $S$, so $S$ is a sequence  converging to $0$ in $H$.  

Let us show that $S$ is qc-dense in $H$. Let $\chi\in\widehat{H}$ and $\chi \ne 0$.  We  need to prove that  $\chi(S)\setminus \T_+\not=\emptyset$. Being a continuous homomorphic image of the compact totally disconnected group $H$, $\chi(H)$ is a  closed totally disconnected subgroup of $\T$. Therefore, $\chi(H)$ must be finite. Hence $\ker \chi$ is an open subgroup of $H$, and consequently it contains a subgroup $W_n$ for some $n\in \N$.  Without loss of generality we will assume that 
\begin{equation}
\label{definition:of:n}
n=\min\{m\in\N:W_m\subseteq\ker\chi\}. 
\end{equation}
Since $\ker\chi\not=H=W_0$  by our assumption, we have $n\ge 1$, and so $n-1\in\N$. 

\medskip
\begin{claim}
$\chi(k_{n-1}v)\not=0$.
\end{claim}
\begin{proof}
Assume the contrary. Then $\chi\restriction_{\hull{k_{n-1}v}} =0$. Since $\hull{v}$ is dense in $H$ and $W_{n-1}$ is an open subset of $H$, 
it follows that $\hull{v}\cap W_{n-1}=\hull{k_{n-1}v}$ is dense in $W_{n-1}$. Now from $\chi\restriction_{\hull{k_{n-1}v}} =0$ and continuity of $\chi$ we conclude that $\chi\restriction_{W_{n-1}}=0$. This gives $W_{n-1}\subseteq \ker\chi$, in contradiction with \eqref{definition:of:n}.
\end{proof}

Since $k_nv\in W_n\subseteq \ker \chi$ by \eqref{definition:of:Wn}
and \eqref{definition:of:n}, we have $k_n\chi(v)=\chi(k_n v)=0$. 
That is, $\hull{\chi(v)}$ is a cyclic group of order at most $k_n$. Since  $\chi(k_{n-1} v)=k_{n-1}\chi(v)\in \hull{\chi(v)}$, the order of  the element $\chi(k_{n-1} v)$ of $\T$ is also at most $k_n$. Since $\chi(k_{n-1}v)\not=0$ by  claim, we can choose an integer $m\le k_n$ such that $\chi(mk_{n-1}v)=m \chi(k_{n-1}v)\not \in\T_+$.  From 
\eqref{definition:of:S}
we conclude that
$mk_{n-1}v\in S$, and so $\chi(S)\setminus \T_+\not=\emptyset$.
\end{proof}

An explicit qc-dense sequence in $\K$ converging to $0$ can be found in \cite[Lemma 4.7]{DS}.
However, that sequence is not contained in $\K_a$. In our next lemma we produce a qc-dense sequence converging to $0$ {\em inside\/} $\K_a$.

\begin{lemma}\label{lemma3} $\K_a$ contains a sequence converging to $0$ that is qc-dense in $\K$.
\end{lemma}

\begin{proof}  We continue using notations from Lemma \ref{lemma2}. Let
$K=\R \times H$  and $u=(1,v)\in K$. Then the cyclic subgroup  $\langle u\rangle$  of  $K$  is discrete and the quotient group $C=K/\langle u\rangle$ is isomorphic to $\widehat{\Q}$  \cite[\S 2.1]{DM}. Therefore, it suffices to prove that $C_a$ contains a sequence converging to $0$ that is qc-dense in $C$.

Let $\pi: K \to C= K/\langle u\rangle$ be the quotient homomorphism. Since $\pi$ is a local homeomorphism, every continuous map $f:[0,1]\to C$ with $f(0)=0_C$ can be lifted to a continuous map $\widetilde f:[0,1]\to K$ with $\widetilde f(0)=0_K$ and $\pi\circ \widetilde f=f$. (A more general statement can be found in \cite[Lemma 1]{Rickert}.)
Therefore, $C_a=\pi(K_a)$.  Since $H$ is zero-dimensional and $\R \times \{0\}$ is arcwise connected, one has $K_a=\R \times \{0\}$, and so $C_a=\pi(\R \times \{0\})$. 

Define $N=\pi(\{0\}\times H)$, and let $f: C \to C/N$ be the quotient homomorphism. By Lemma \ref{lemma2}  and Fact \ref{continuous:images:are:qc-dense},
there exists a converging to 0 sequence $S'$ in the subgroup $\hull{\pi(0,v)}$ of $N$ such that $S'$ is qc-dense in $N$.  As 
$$
\pi(0,v)=\pi((-1,0)+(1,v))=
\pi(-1,0)+\pi(u)=
-\pi(1,0)\in \pi(\R \times \{0\})= C_a,
$$ 
one has $S' \subseteq \pi (\hull{(0,v)})\subseteq  C_a$. 

 With $T=\left\{\frac{1}{2n}: n\in \N, n\ge 1\right\}$  define $S''=\pi(T\times \{0\})\subseteq \pi(\R \times \{0\})=C_a$. Clearly, $S''$ is a sequence converging to $0$. Since 
$C/N\cong K/(\Z\times H) \cong \T$ and the  composed isomorphism $C/N\to \T$ sends $f(S'')$ to $\varphi(T)$, from Example \ref{Example1} we conclude that $f(S'')$  is qc-dense in $C/N$.

Since $S'$ and $S''$ are sequences converging to $0$ in $C$, so is $X= S'\cup S''$. By our construction, $X\subseteq C_a$. So it remains only to prove that
$X$ is qc-dense in $C$. Suppose that $\chi\in\widehat{C}$ and $\chi(X)\subseteq \T_+$. Since $\chi\restriction_N\in \widehat{N}$, 
$\chi\restriction_N(S')=\chi(S')\subseteq \chi(X)\subseteq \T_+$ and $S'$ is qc-dense in $N$, we have $\chi\restriction_N=0$. Therefore, $\chi=\xi \circ f$ for some
$\xi\in\widehat{C/N}$. In particular, $\xi(f(S''))\subseteq \xi(f(X))=\chi(X)\subseteq \T_+$. Since $f(S'')$ is qc-dense  in 
$C/N$, it follows that $\xi=0$. This gives $\chi =0$. Therefore, $X$ is qc-dense in $C$. 
\end{proof}

\section{Proof of Theorems  \ref{super-sequence:theorem} and \ref{MainThm}}\label{Sec2}\label{proofs}

The following definition is an adaptation to the abelian case of \cite[Definition 4.5]{DS1}:

\begin{definition}
\label{definition:of:fan}
Let $\{G_i:i\in I\}$ be a family of abelian topological groups. For every $i\in I$ let $X_i$ be a subset of $G_i$. Identifying each $G_i$ with a subgroup of the direct product 
$G=\prod_{i\in I} G_i$ in the obvious way, define $X=\bigcup_{i\in I} X_i\cup\{0\}$, where  $0$ is the zero element of $H$. We will call $X$ 
the {\em fan\/} of the family $\{X_i:i\in I\}$ and will denote it by $\mathrm{fan}_{i\in I}(X_i, G_i)$. 
\end{definition}

The proof of the following lemma is straightforward.

\begin{lemma}\label{lemma1}
{\rm \cite[Lemmas 4.3 and 4.4]{DS}}
 Let $\{G_i:i\in I\}$ be a family of abelian topological groups, and let $G=\prod_{i\in I} G_i$. For every $i\in I$ let $X_i$ be a subset of $G_i$, and let $X=\mathrm{fan}_{i\in I}(X_i,G_i)$. Then: 
\begin{itemize} 
\item[(i)]  if $X_i$ is a sequence converging to $0$ in $G_i$, then $X$ is a super-sequence in $G$ converging to $0$. 
\item[(ii)] if $X_i$ is a qc-dense subset of $G_i$ for each $i\in I$, then $X$ is qc-dense in $G$. 
\end{itemize}
\end{lemma}

\begin{lemma}
\label{corollary1} 
For every cardinal $\kappa$ there exists a super-sequence $S\subseteq (\K^\kappa)_a$ converging to  $0$ that is qc-dense in $\K^\kappa$.
\end{lemma}

\begin{proof}
Write $\K^\kappa$ as $\K^\kappa= \prod_{\alpha<\kappa}G_\alpha$, where $G_\alpha$ is the $\alpha$'s copy of $\K$. By Lemma \ref{lemma3}, for every $\alpha\in\kappa$ there is a sequence $S_\alpha $ in $(G_\alpha)_a$ converging  to 0  that is qc-dense in $G_\alpha$. By Lemma \ref{lemma1}(i), $S=\mathrm{fan}_{\alpha\in \kappa}(X_\alpha,G_\alpha)$ is a super-sequence in $\K^\kappa$ converging to 0. By Lemma \ref{lemma1}(ii), $S$ is qc-dense in $\K^\kappa.$  Finally, note that $S\subseteq \bigoplus_{\alpha\in\kappa} (G_\alpha)_a\subseteq (\K^\kappa)_a$.  
\end{proof}

\begin{proof}[\bf Proof of Theorem  \ref{super-sequence:theorem}]
(i)$\to$(ii) Let $S\subseteq G_a$ be a super-sequence that is qc-dense in $G$. Then
the subgroup $H$ of $G$ generated by $S$ is dense in $G$ by Fact \ref{qc-dense:implies:topological:generation}. Since $H\subseteq G_a$, it follows that $G_a$ is dense in $G$ as well. Thus $G$ is connected.

(ii)$\to$(i)
There exits a continuous surjective homomorphism  $\pi: \K^\kappa \to G$ for some cardinal 
$\kappa$ (see, for example, the proof of \cite[Theorem 3.3]{DS}). Let $S$ be as in the conclusion of Lemma \ref{corollary1}. 
Since $S$ is qc-dense in $\K^\kappa$, $\pi(S)$ is qc-dense in $G$ by Fact \ref{continuous:images:are:qc-dense}. 
Since a finite set cannot be qc-dense in an infinite compact group (\cite{Diss}; this also follows from Theorem \ref{weight:of:qc-dense:set}),
$\pi(S)$ must be infinite. Since $S$ is a super-sequence converging to 0 such that $\pi(S)$ is infinite, 
$\pi(S)$ must be a super-sequence converging to 0  by \cite[Fact 4.3]{DS1} (see \cite[Fact 12]{Sh} for the proof).
Finally note that $\pi(S)\subseteq \pi(\K^\kappa_a)\subseteq G_a$.
\end{proof}

\begin{proof}[\bf Proof of Theorem  \ref{MainThm}] 
We  have $G = \R^n \times K$, where $K$ is a compact connected group \cite{DPS}.  Since $\R^n$ is arcwise connected, one has $G_a= \R^n \times K_a$. 
From Theorem \ref{super-sequence:theorem} and Fact \ref{connection:between:determination:and:qc-density} we conclude that $K_a$ determines $K$. Hence $G_a= \R^n \times K_a$ determines $G= \R^n \times K$.
\end{proof}

\section{Proof of Theorem \ref{topologically:generating:sequence}}

In the sequel we denote by $H'$ the commutator subgroup of a group $H$. 

Our next lemma shows that qc-density can be essentially studied in the abelian context. 
\begin{lemma}
\label{non-abelian:qc-dense} 
Let $H$ be a topological group, and let $G$ denote the quotient $H/\overline{H'}$, where $\overline{H'}$ is the closure of
$H'$ in $H$.  Let $\pi: H\to G$ denote the canonical map. Then a subset $E$ of $H$ is qc-dense in $H$ if and only if $\pi(E)$ is qc-dense in $G$. 
\end{lemma}

\begin{proof}
The ``only if'' part follows from Lemma \ref{continuous:images:are:qc-dense}. To prove the ``if'' part, assume that  $\pi(E)$ is qc-dense in $G$,  and let $\chi:H \to \T$ be a continuous homomorphism such that $\chi(E) \subseteq \T_+$. Since $\T$ is  abelian and Hausdorff, $\overline{H'}\subseteq \ker \chi$, so $\chi= \xi\circ \pi$ for some  character  $\xi: G \to \T$. Since $\xi(\pi(E))= \chi(E)\subseteq \T_+$ and $\pi(E)$ is qc-dense in $G$, we conclude that $\xi$ is trivial, and so $\chi$ is trivial too. This proves that $E$ is qc-dense in $H$. 
\end{proof}

Recall that a group $L$ is called {\em perfect\/} if $L'=L$. 

\begin{corollary}\label{corollary:product} Let $L$ be a perfect topological group, $G$ an abelian topological group and $H = G \times L$. Then a subset 
$E$ of  $G$ is qc-dense in $G$ if and only if the subset $E \times \{e_L\}$ of the group $H$ is qc-dense in $H$.
\end{corollary}

\begin{proof}
Since $H'=\{0_G\} \times L$, we have $H'=\overline{H'}$ and $G\cong H/H'= H/\overline{H'}$. Now the conclusion of our corollary
follows from Lemma \ref{non-abelian:qc-dense} applied to the projection $\pi: H \to G$.
\end{proof}

\begin{example}
\label{Example:2}
Let $T=\left\{\frac{1}{2n}: n\in \N, n\ge 1\right\}\cup\{0\}$. Assume that $L$ is a finite simple non-abelian group and $G = \T \times L$.
  Then  {\em $S = \varphi(T)\times \{e_L\}$ is a qc-dense sequence in $G$ converging to $e_G$\/}. Indeed, since  all simple non-abelian groups are perfect, 
this follows from Example \ref{Example1} and Corollary \ref{corollary:product}.
\end{example}

\begin{proof}[\bf Proof of Theorem \ref{topologically:generating:sequence}] 
(i)$\to$(iii) Let $S\subseteq G_a$ be a suitable set for $G$. Then the subgroup $H$ of $G$ generated by $S$ is dense in $G$. Since $H\subseteq G_a$, it follows that $G_a$ is dense in $G$ as well. Thus $G$ is connected.

(iii)$\to$(ii)
According to \cite[Theorem 3.3]{DS1}, there exist a cardinal $\kappa$ and a continuous surjective group homomorphism $\pi: \K^\kappa\times L\to G$, where $L$ is a direct product of simple Lie groups. In particular, $L$ is perfect.  By Lemma \ref{corollary1}, there exists a super-sequence $E$ in $(\K^\kappa)_a$ converging to $0$ that is qc-dense in $\K^\kappa$.  By Corollary \ref{corollary:product},  $S'=E\times  \{e_L\}$ is a converging to $(0,e_L)$ super-sequence in $(\K^\kappa)_a\times \{e_L\}$ that is qc-dense in $\K^\kappa\times L$. Since $S'$  is qc-dense in $\K^\kappa \times\{e\}$, it topologically generates $\K^\kappa\times\{e\}$ by Fact \ref{qc-dense:implies:topological:generation}. By Theorem \ref{Hofmann:Morris} and Remark \ref{suitable:remark},
there exists a super-sequence $S''\subseteq \{0\}\times L$  converging to $(0,e_L)$ that topologically generates $\{0\}\times L$.  Then $S=S'\cup S''$ is a  super-sequence converging 
to $(0,e_L)$ that is qc-dense in $\K^\kappa\times L$ and topologically generates $\K^\kappa\times L$. 
Therefore,  $X=\pi(S)$ is a super-sequence (by \cite[Fact 4.3]{DS1}; see \cite[Fact 12]{Sh} for the proof)
that is qc-dense in $G$ (by Fact \ref {continuous:images:are:qc-dense}) and topologically generates $G$. 
Since $L$ is arcwise connected and $S\subseteq (\K^\kappa)_a\times L$,  we have 
$X=\pi(S)\subseteq \pi((\K^\kappa)_a\times L)\subseteq G_a$. Applying Remark \ref{suitable:remark}, we conclude that $X\setminus\{e\}$ is a suitable set for $G$. Since $X$ is qc-dense in $G$, so is $X\setminus\{e\}$. 
If $X$ is infinite, we are done. 

Assume now that $X$ is finite. Since $S'$ is qc-dense in $\K^\kappa\times\{e\}$, $\pi(S')$ is qc-dense in $K=\pi(\K^\kappa\times\{e\})$ by Fact \ref {continuous:images:are:qc-dense}. Since $\pi(S')\subseteq \pi(S)=X$, the set $\pi(S')$ is finite. From Theorem \ref{weight:of:qc-dense:set} we conclude that the (compact abelian) group $K$ must finite as well. Being a continuous image of the connected group $\K^\kappa\times\{e\}$, the group $K$ is connected. Hence,  $K$ is trivial, and so
$G = \pi(\K^\kappa\times L)= \pi(\{0\}\times L) \subseteq G_a\subseteq G$  because $L$ is pathwise connected. This yields $G=G_a$.

Let $f:[0,1]\to G_a$ be a continuous map such that 
$f(0)=e_G$ 
and $f(1)\not=e_G$.
Define $c=\inf\{t\in [0,1]: f(t)\not=e_G\}$. Then $f(c)=e_G$ by continuity of $f$ and the choice of $c$. For every $n\in\N$ choose $c_n\in[0,1]$ such that $c<c_n<c+1/n$ and $f(c_n)\not=e_G$. Since $\{c_n:n\in\N\}$ converges to $c$ and $f$ is continuous, the sequence $X_0=\{f(c_n):n\in\N\}$ converges to $f(c)=e_G$.
Since $f(c_n)\not=e_G$ for every $n\in\N$, we conclude that 
$X_0$ is infinite.
Now $X'=X \cup X_0\cup \{e_G\}$
is an infinite super-sequence converging to $e_G$, and so
$X'\setminus\{e_G\}\subseteq G_a$ is an infinite suitable set for $G$ by Remark \ref{suitable:remark}. Since 
$X$ is qc-dense in $G$ and topologically generates $G$, $X'$ has the same properties.

The implication (ii)$\to$(i) is trivial.
\end{proof}

\end{document}